\documentclass{dcds-b}
\usepackage{a4}
\usepackage{paralist}

\def\currentvolume{XX}
 \def\currentissue{X}
  \def\currentyear{2011}
   \def\currentmonth{Month X}
    \def\ppages{1--x}
     \def\DOI{yy.yyyy/dcdsb.2010.yy.y}

\usepackage{amsmath, amssymb,color,ulem}

\newtheorem{theorem}{Theorem}

\newtheorem{proposition}[theorem]{Proposition}

\newcommand{\ud}{\mathrm{d}}

\newcommand{\A}{\mathcal{A}}
\newcommand{\B}{\mathcal{B}}

\newcommand{\T}{\mathcal{T}}

\newcommand{\ra}{\rightarrow}

\renewcommand{\theequation}{\thesection.\arabic{equation}}
\newcommand{\eqdef}{{\ \stackrel{\mathrm{def}}{=}\ }}

\title[Steady states in hierarchical populations]{Steady states in hierarchical
structured populations with distributed states at birth}
\author[J.\,Z. Farkas, P. Hinow]{}
\date{\today}

\email{jzf@maths.stir.ac.uk}
\email{hinow@uwm.edu}

\begin{document}
\maketitle

\centerline{\scshape J\'{o}zsef Z. Farkas }
{\footnotesize
 \centerline{Institute of Computing Science and Mathematics}
 \centerline{University of Stirling}
\centerline{Stirling, FK9 4LA, United Kingdom}  
} 
\medskip
\centerline{\scshape Peter Hinow }
{\footnotesize
 \centerline{Department of Mathematical Sciences}
 \centerline{University of Wisconsin -- Milwaukee}
 \centerline{P.O. Box 413, Milwaukee, WI 53201-0413, USA}
} 
\bigskip
\begin{abstract}
\begin{sloppypar}
We investigate steady states of a
quasilinear first order hyperbolic partial integro-differential equation. The
model describes the evolution of a hierarchical structured population with
distributed states at birth. Hierarchical size-structured models describe the
dynamics of populations when individuals experience size-specific environment.
This is the case for example in a population where in\-di\-vi\-du\-als exhibit
cannibalistic behavior and the chance to become prey (or to attack) depends on
the individual's size. The other distinctive feature of the model is that
individuals are recruited into the population at arbitrary size. This amounts to
an infinite  rank integral operator describing the recruitment process. First we
establish conditions for the existence of a positive steady state of the model.
Our method uses a fixed point result of nonlinear maps in conical shells of
Banach spaces.  Then we study stability properties of steady states for the
special case of a separable growth rate using results from the theory of
positive operators on Banach lattices. 
\end{sloppypar}
\end{abstract}
\maketitle

\section{Introduction}\label{section:introduction}
Classic population models often assume that individuals experience
\textit{scramble} competition. This means that all individuals in the
population have equal chances in the competition for resources such as food,
light, space etc., see e.g.~\cite{CUS,GuMC,MD,P,WEB}. In many species, however,
competition for resources that determine individual mortality and fertility is
based on some hierarchy in the population which is related to individuals size
or to any other variable that characterizes physiological structure. Significant
amount of interest has been paid to understand the dynamics of populations that
exhibit \textit{contest} competition. Structured population models are a useful
tool to study intra-specific contest competition. Both (time-)discrete (see
e.g.~\cite{JC,JC2}) and continuous (see e.g.~\cite{ADH,AI,CS3,CS2,C1,FH2,FH3})
models have been formulated and analysed to this end. Most of the nonlinear
models in the literature incorporate environmental feedback through some form
of density dependence in the vital rates. Our goal in this paper is to carry out
a qualitative analysis of a continuous model which incorporates quite general
nonlinearities. 

We study the following quasilinear partial integro-differential equation
\begin{align}
\frac{\partial}{\partial t}p(s,t)+ \frac{\partial}{\partial
s}\left(\gamma(s,P(t))p(s,t)\right) 
& =-\mu(s,E(s,p))p(s,t)  \nonumber \\
&\phantom{=} +\int_0^m\beta(s,y,E(y,p))p(y,t)\,\ud y,\label{equation1} \\
\gamma(0,P(t))p(0,t)&=0,\label{equation2} \\
E(s,p)  =\alpha\int_0^s w(r)p(r,t)\,\ud r&+\int_s^mw(r)p(r,t)\, \ud
r,\label{equation3} \\
P(t) & =\int_0^m\kappa(s)p(s,t)\,\ud s,\label{equation4}
\end{align}
with the initial condition
\begin{equation*}
p(s,0)=p_0(s).
\end{equation*}
This model describes the long term dynamics of a population of a sufficiently
large size living in a closed
habitat. The function $p=p(s,t)$ denotes the density of individuals of size (or
other developmental
stage) $s$ at time $t$. It is assumed that individuals may have different sizes
at birth and therefore $\beta(s,y,\,\cdot\,)$ denotes the rate at which 
individuals of size $y$ ``produce'' individuals of size $s$. Hence the non-local
integral term in Equation \eqref{equation1} represents the recruitment of
individuals into the population. $\gamma$ denotes the size-specific growth rate,
while $\mu$ stands for the mortality rate. We assume that individual growth is
also regulated by a weighted population size $P$, for example due to
competition. Mortality however, depends on the size-specific environment $E$,
for example due to cannibalism. The parameter $\alpha$ is related to
the strength of the hierarchy in the population. We note that if $\alpha=1$ our
model reduces to the scramble competition model that in the special case 
$\kappa\equiv w\equiv1$ was considered in \cite{FGH}. We make the following
general assumptions on the model ingredients
\begin{equation}\label{assumptions}
\begin{aligned}
& \mu\in C^1([0,m]\times [0,\infty)),\quad \beta\in C^1([0,m]\times
[0,m]\times [0,\infty)), \\
&\gamma\in C^2([0,m]\times [0,\infty)), \quad w,\,\kappa\in
L^\infty(0,m),\\
&\beta,\,\alpha,\,w,\,\mu\geq 0, \quad \gamma,\,\kappa>0.
\end{aligned}
\end{equation}
Notice that we make no requirement
that $\beta(s,y,E(y,p))=0$ if $y<s$ although this seems natural from a
biological point of view.
In the remarkable paper \cite{CS2}, Calsina and Salda\~{n}a studied
well-posedness of a very general size-structured model with distributed states
at birth. They established global existence and uniqueness of solutions
using results from the theory of nonlinear evolution equations. 
Model \eqref{equation1}-\eqref{equation4} is a special case of the general
model treated in \cite{CS2}, however, in \cite{CS2} qualitative questions were
not addressed.  In contrast to \cite{CS2}, our paper focuses on the existence
and local asymptotic stability of equilibrium solutions of system
\eqref{equation1}-\eqref{equation4} with particular regards to the effects 
of  the distributed states at birth (previously we addressed simpler models 
without hierarchical structure in
\cite{FGH,FH2,FH3}). Earlier, in \cite{HC} Henson and Cushing studied
continuous age-structured hierarchical models. They compared models with
scramble versus contest competition for a limited resource. 
In particular the equilibrium levels for the two modes of competition were
analyzed. Crucially however, the models in \cite{HC} incorporate vital rates
that do not depend explicitly on the structuring variable.  

In Section \ref{section:Existence}, we will establish conditions for the
existence of
positive steady states of our model. The question of the existence of
non-trivial steady states is difficult mainly for two reasons. Firstly,
due to hierarchy in the population related to individual size, individual
mortality and fertility depend on the size specific environment $E$. This
environmental feedback yields an infinite dimensional nonlinearity in the model
equations, in contrast to Gurtin-MacCamy type models \cite{GuMC}, where the
vital rates depend on a weighted total population size, or on a finite number of
such variables. Secondly, as individuals may be recruited into the population
at all possible sizes, a recruitment operator of infinite  rank 
arises. This means that  the steady state equation cannot be solved
explicitly. In \cite{FGH} we overcame this issue for a simpler model where the
model ingredients only depended on the total population size by using results
from the spectral theory of positive operators. Unfortunately this approach
cannot be extended to the model considered here. Therefore we devise a
different approach, based on fixed point results for nonlinear maps in conical
shells of a Banach space, see \cite{Amann,Deimling}. This approach was used
before to treat age-structured models (also with diffusion), see \cite{P,W2},
where every
individual enters the population at a single state, namely at age zero. The
method is based on the construction of an appropriate nonlinear map that
requires the 
implicit solution of the steady state equation. However, the solution of the
steady state equation of our model is not available. Therefore we need to
construct an appropriate sequence of recruitment processes of finite rank for
which we can solve the corresponding steady state problems. Then we show that
the steady states constructed from the fixed points of the sequence of the
nonlinear maps, have actually a convergent subsequence and we show that the
limit point is actually a steady state of the original problem
\eqref{equation1}-\eqref{equation4}.

In Section \ref{section:Asymptotic} we focus on the asymptotic behavior of solutions of the
model. A positive quasicontraction semigroup describes the evolution of
solutions of the system linearized at an equilibrium solution. We establish a
regularity property of the governing linear semigroup in Proposition \ref{reg}
that allows in principle to address stability questions of positive
equilibrium solutions of \eqref{equation1}-\eqref{equation4}. However, even
the point spectrum of the linearized  semigroup generator cannot be
characterized explicitly via zeros of an associated cha\-rac\-te\-ris\-tic
function. 
This is because the eigenvalue equation cannot be solved explicitly due to the
infinite dimensional nonlinearity in the original model  and the very general
recruitment process. We will overcome this issue by using compact positive
perturbations of the semigroup generator and rank one perturbations of the 
general recruitment term. This allows us to arrive at stability/instability
conditions for the steady states of our model.

\section{Existence of positive equilibrium solutions}\label{section:Existence}

In this section we will discuss the existence of steady states of model
\eqref{equation1}-\eqref{equation4}.  
We define the nonlinear operator 
\begin{equation*}
\Psi\, :\, W^{1,1}(0,m)\to L^1(0,m)
\end{equation*}
by
\begin{equation}\label{eq1}
\Psi(q)= \frac{\partial}{\partial
s}\left(\gamma(s,Q)q(s)\right)+\mu(s,E(s,
q))q(s)-\int_0^m\beta(s , y , E(y,q))q(y)\, \ud y,
\end{equation}
where 
\begin{equation*}
Q=\int_0^m\kappa(s)q(s)\,\ud s.
\end{equation*}
It is clear that $p_*\in W^{1,1}(0,m)$ is a steady state of
\eqref{equation1}-\eqref{equation4}
if and only if $\Psi(p_*)=0$ and $p_*(0)=0$.

Our aim here is to apply a fixed point result, see e.g.~Theorem 12.3
in \cite{Amann} or Theorem A in \cite{P}. Its proof uses the
Leray-Schauder degree theory for compact perturbations of the identity in
infinite dimensional Banach spaces. 

\begin{theorem}\label{Theorem:FixedPoint}
Let $(\mathcal{X},||\,\cdot\,||_{\mathcal{X}})$ be a Banach space,
$\mathcal{K}\subset\mathcal{X}$ a closed convex cone and
$\mathcal{K}_r=\mathcal{K}\cap B_r(0)$, 
where $B_r(0)$ denotes the ball of radius $r$ centered at the origin. Let
$\Phi:\, \mathcal{K}_r\to \mathcal{K}$ continuous such that
$\Phi(\mathcal{K}_r)$ 
is relatively compact. Assume that
\begin{enumerate}
\item\label{FPc1} $\Phi x\ne \lambda x\quad \text{for all}\quad
||x||_{\mathcal{X}}=r,\,\lambda>1.$
\item\label{FPc2} There exists a $\rho\in (0,r)$ and
$k\in\mathcal{K}\setminus\{0\}$ such that 
\begin{equation*}
x-\Phi x\ne\lambda k\quad\text{for all}\quad
||x||_{\mathcal{X}}=\rho,\,\lambda>0.
\end{equation*}
\end{enumerate}
Then $\Phi$ has at least one fixed point $x_*$ in the shell 
\begin{equation*}
\mathcal{S}_{\rho,r}=\{x\in\mathcal{K}\,:\,\rho\le||x||_{\mathcal{X}}\le r\}.
\end{equation*}
\end{theorem}

\subsection{The case of a separable fertility
function}\label{subsection:Existence1}
First we show in Theorem \ref{theorem:ex-sep} that for a finite rank fertility
process the problem admits a positive steady state under biologically
meaningful conditions on the model ingredients. Later we treat the general case
in Theorem \ref{theorem:existence_full}.
\begin{theorem}\label{theorem:ex-sep}
Assume that 
\begin{equation}\label{specbeta}
\beta(s,y,E(y))=\sum_{j=1}^l\beta_j(s)\bar{\beta}_j(y,E(y))
\end{equation}
with continuous functions $\beta,\,\bar{\beta}$ and there
exists a $j\in\{1,...,l\}$ such that 
\begin{equation}\label{sscond1}
\int_0^m\bar{\beta}_j(s,0)F_j(s,\mathbf{0},0)\,\ud s>1,
\end{equation}
where $\mathbf{0}$ is the zero function and
\begin{equation*}
F_j(s,E(s),P)=\int_0^s\exp\left\{-\int_x^s\frac{\mu(r,E(r))+\gamma_s(r,P)}{
\gamma(r,P)}\,\ud r\right\}\frac{\beta_j(x)}{\gamma(x,P)}\,\ud x.
\end{equation*}
Let $F$ be a bounded and measurable function that satisfies
\begin{equation*}
F(s,H(s),P) \ge F_j(s,H(s),P)
\end{equation*}
for all $j$, $s\in [0,m],\,H\in L_+^1(0,m)$ and $P>0$, and $c$ be a constant
such that
\begin{equation}\label{cinequality}
\kappa(s)\ge c\sum_{k=1}^l\bar{\beta}_k(s,H(s)).
\end{equation}
Suppose that there exists an $R>0$ such that for all
$(H,P)\in L_+^1(0,m)\times \mathbb{R}_+$  with $||H||_{L^1}+P>R$
we have 
\begin{equation}\label{sscond3}
\int_0^m\kappa(s)F(s,H(s),P)\,\ud s\le c.
\end{equation}
Then the model \eqref{equation1}-\eqref{equation4} admits a positive steady
state $p_*$.
\end{theorem}
\begin{sloppypar}
\noindent {\bf Proof.} Let $\mathcal{X}=L^1(0,m)\oplus \,l^1$ with norm
$||\,\cdot\,||_{\mathcal{X}}=||\,\cdot\,||_{L^1}+||\,\cdot\,||_{l^1}$. We will
use the notation \mbox{$x=(H,\mathbf{P})$}, where $H\in L^1(0,m)$ and
\mbox{$\mathbf{P}\in l^1$} or \mbox{$x=(H,P^0,\mathbf{P}')$} where $P^0$ is just
the first component of $\mathbf{P}$ and \mbox{$\mathbf{P}'=(P^1,P^2,\dots)$}. We
consider elements of $\mathbb{R}^n$ to be in $l^1$ by the trivial embedding. We
denote by $\mathcal{K}=\left(L^1(0,m)\oplus l^1\right)_+$ the positive cone of
$\mathcal{X}$ which is closed and convex and denote
$\mathcal{K}_r=\mathcal{K}\cap B_r(0)$, where $r$ has yet to be chosen. Without
loss of generality we may assume that the indices have been assigned such that
condition \eqref{sscond1} holds for $j=1$. 
\end{sloppypar}

We note that for a fertility function $\beta$ of the form \eqref{specbeta} the 
non-trivial time independent solution of model
\eqref{equation1}-\eqref{equation4} can be found (if it exists) as
\begin{equation}\label{sssol1}
p_*(s)=\sum_{j=1}^lP_*^jF_j(s,E_*(s),P^0_*),
\end{equation}
where
\begin{equation}\label{sssol4}
\begin{aligned}
P_*^j & = \int_0^m\bar{\beta}_j(s,E_*(s))p_*(s)\,\ud s,\quad
j\in\{1,\dots,l\}, \\
E_*(s) & = \,\alpha\int_0^sw(r)p_*(r)\,\ud r+\int_s^mw(r)p_*(r)\,\ud r, \\
P^0_* & =\int_0^m\kappa(s)p_*(s)\,\ud s.
\end{aligned}
\end{equation}
We define a nonlinear map $\Phi^l\,:\,\mathcal{K}_r\to\mathcal{K}$ by
\mbox{$\Phi^l=\left(\Phi_1^l,\Phi_2^l\right)$}, where 
\begin{align*}
&\Phi_1^l(H,\mathbf{P})(s) \\ 
&=\sum_{j=1}^l
P^j\left(\alpha\int_0^sw(r)F_j(r,H(r),P^0)\,\ud
r+\int_s^mw(r)F_j(r,H(r),P^0)\,\ud
r\right),  \\
&\Phi_2^l(H,\mathbf{P})=
\left(\sum_{j=1}^lP^j\int_0^m\kappa(s)F_j(s,H(s),P^0)\,\ud
s,\right. \\
& \phantom{\Phi_2^l(H,\mathbf{P})=(} 
\sum_{j=1}^lP^j\int_0^m\bar{\beta}_1(s,H(s))
F_j(s,H(s),P^0)\,\ud s,\,\dots\,, \\
&\phantom{\Phi_2^l(H,\mathbf{P})=(}\left.\sum_{j=1}^lP^j\int_0^m\bar{\beta}_l(s,
H(s)
)F_j(s,H(s),P^0)\,\ud s, 0, \dots  \right). 
\end{align*}
\begin{sloppypar}
Although at this moment we consider only a single fertility
function $\beta$
(namely the finite sum of separable functions from equation \eqref{specbeta}),
we use the superscript notation $\Phi^l$ in anticipation of a sequence of such
maps that we shall employ in the proof of Theorem \ref{theorem:existence_full}
below. Note that \mbox{$\Phi_1^l(H,\mathbf{P})\ge 0$} and since
\begin{equation*}
(\Phi_1^l)'(H,\mathbf{P})(s) =\sum_{j=1}^l
P^j(\alpha-1) w(s)F_j(s,H(s),P^0), 
\end{equation*}
also $\Phi_1^l(H,\mathbf{P})\in W^{1,1}(0,m)$. Likewise 
\mbox{$\Phi_2^l(H,\mathbf{P})\ge 0$} and the
nonzero entries of $\Phi_2^l(H,\mathbf{P})$ are in
$\mathbb{R}^{l+1}_+$. Thus $\Phi$ maps indeed into $\mathcal{K}$ and has a
relatively compact image. It
can be seen that $(E_*,\mathbf{P}_*)\in\mathcal{K}_r$ is a fixed point of the
nonlinear map $\Phi^l$ if and only if the function $p_*$ defined via 
equations \eqref{sssol1}-\eqref{sssol4} is a steady state of problem
\eqref{equation1}-\eqref{equation4}.
\end{sloppypar}

Assume now that for $(H,P^0,\mathbf{P}')\in\mathcal{K}_r$ such that
\begin{equation*}
||(H,P^0,\mathbf{P}')||_{\mathcal{X}}=||H||_{L^1}+P^0+||\mathbf{P}'||_{l^1}
=r
\end{equation*} 
and for some $\lambda>1$ we have 
\begin{equation*}
\Phi^l(H,P^0,\mathbf{P}')=\lambda (H,P^0,\mathbf{P}'),
\end{equation*} 
that is 
\begin{equation}\label{cond1}
\Phi^l_1(H,P^0,\mathbf{P}')=\lambda H,\quad
\Phi^l_2(H,P^0,\mathbf{P}')=\lambda(P^0,\mathbf{P}').
\end{equation}
The second equation in \eqref{cond1} can be written as
\begin{align}
\lambda P^0 & =\sum_{j=1}^lP^j\int_0^m
\kappa(s)F_j(s,H(s),P^0)\,\ud s, \label{firsteq}\\
\lambda P^1 & = \sum_{j=1}^lP^j\int_0^m\bar{\beta}_1(s,H(s))
F_j(s,H(s),P^0)\,\ud s, \label{secondeq}\\
\dots & \nonumber \\
\lambda P^l & = \sum_{j=1}^lP^j\int_0^m\bar{\beta}_l(s,H(s))
F_j(s,H(s),P^0)\,\ud s. \label{leq}
\end{align}
It follows from equations \eqref{cond1}-\eqref{leq} that we may assume that
$P^j\ne 0$ for $j=0,1,2,\dots,l$. 

From equations \eqref{secondeq}-\eqref{leq} we obtain
\begin{equation}\label{l1norm}
\lambda||\mathbf{P}'||_{l^1}=\sum_{k=1}^l\sum_{j=1}^lP^j\int_0^m\bar{\beta}_k(s,
H(s))F_j(s,H(s),P^0)\,\ud s.
\end{equation}
Combining equations \eqref{firsteq} and \eqref{l1norm} we obtain
\begin{equation*}
||\mathbf{P}'||_{l^1}=P^0\left(\frac{\displaystyle\sum_{j=1}
^lP_j\int_0^m\displaystyle\sum_{k=1}^l\bar{\beta}_k(s,H(s))F_j(s,
H(s),P^0)\,\ud s}
{\displaystyle\sum_{j=1}^lP_j\int_0^m\kappa(s)F_j(s,H(s),P^0)\,\ud
s}\right).
\end{equation*}
This, combined with inequality
\eqref{cinequality} implies that
\begin{equation}\label{estim2}
c\,||\mathbf{P}'||_{l^1}\le P^0.
\end{equation}
However, by choosing $r=||H||_{L^1}+P^0+||\mathbf{P}'||_{l^1}
\ge ||H||_{L^1}+P^0>R$ sufficiently large and using condition \eqref{sscond3}, 
we have from equation
\eqref{firsteq}
\begin{align*}
P^0 & <\sum_{j=1}^lP^j\int_0^m\kappa(s)F_j(s,H(s),P^0)\,\ud s \\
& \le||\mathbf{P}'||_{l^1}\int_0^m\kappa(s)F(s,H(s),P^0)\,\ud s \le c
||\mathbf{P}'||_{l^1}\le P^0,
\end{align*}
a contradiction. Thus condition \eqref{FPc1} of Theorem
\ref{Theorem:FixedPoint} is established. 

Let us now define $k=(0,(1,\dots,1,0,\dots))\in\mathcal{K}\setminus\{0\}$
with $l+1$ entries $1$ and assume that for some $\lambda>0$ and $\rho>0$ we have
for all $(H,\mathbf{P})$ with $||(H,\mathbf{P})||_{\mathcal{X}}=\rho$  
\begin{equation*}
(H,\mathbf{P})-\Phi^l(H,\mathbf{P})=\lambda\,k,
\end{equation*} 
that is 
\begin{equation*}
H-\Phi^l_1(H,\mathbf{P})=0,\quad \mathbf{P}-\Phi_2^l
(H,\mathbf{P})=(\lambda,\dots,\lambda).
\end{equation*}
The latter equation  can be written as
\begin{equation}\label{cond2}
(\mathbf{I}-\mathbf{B}(H,P^0))\cdot\mathbf{P}=(\lambda,\dots,\lambda) ,
\end{equation}
where the $(l+1)\times(l+1)$ matrix $\mathbf{B}$ has elements $B_{00}=0$ and 
\begin{equation*}
B_{ii}=\int_0^m\bar{\beta}_i(s,H(s))F_i(s,H(s),P^0)\,\ud
s\quad \text{for}\quad i=1,\dots,l.
\end{equation*}
It follows from condition \eqref{sscond1} and the continuity of
$\bar{\beta}_{1}$ and $F_{1}$ that $B_{11}>1$ for all $||(H,P^0)||=\rho$ for
some small enough value $\rho>0$. This renders the left hand side of
\eqref{cond2} negative and yields a contradiction. Condition \eqref{FPc2} of
Theorem \ref{Theorem:FixedPoint} is established and the proof is now completed.
\hfill $\Box$
   
\subsection{General fertility function}\label{subsection:Existence2} 
Every fertility function $\beta$ of the required regularity $C^1$
(in all its arguments, see \eqref{assumptions}) can be
written as a limit of partial sums of separable functions 
\begin{equation}\label{betaseries}
\beta^l(s,y,E(y))=
\sum_{k=1}^{l}\beta^l_k(s)\bar{\beta}^l_k(y,E(y)),
\end{equation}
with
\begin{equation*}
\lim_{l\to\infty}||\beta(\,\cdot\,,\,\cdot\,,E(\,\cdot\,))-\beta^l(\, \cdot\,,
\cdot\,,E(\,\cdot\,))||_ {L^\infty([0,m]^2)}=0
\end{equation*}
for every $E\in L^1(0,m)$. This can be achieved by partitioning the interval
$[0,m]$ into $l$ subintervals by $y_k^l=k\frac{m}{l}$, $k=0,\dots,l$ and
setting 
\begin{equation*}
\beta_{k+1}^{l}(s) = \beta(s, y^l_k, E(y^l_k)) \quad
\textrm{and} \quad
\bar{\beta}_{k+1}^{l}(y,E(y)) = \chi_{[y^l_{k-1},y^l_k]}(y)
\end{equation*}
for  $k=1,\dots,l$ where $\chi_{[y^l_{k-1},y^l_k]}$ denotes the indicator
function
of the respective interval. Then the partial sums converge in the supremum norm
as the number of subintervals increases to infinity and there is a uniform
bound of the derivatives in the $s$-direction
\begin{equation}\label{excond4}
||\beta_k||_{C^1([0,m])}\le M\quad 
\end{equation}
for all $k$, uniformly for arbitrary $E\in L^1$. The latter is
because of the boundedness of all derivatives of $\beta$.

\begin{theorem}\label{theorem:existence_full}
Assume that for fixed model ingredients $\mu,\gamma,\beta$ and $\kappa$ there
exists a bounded and measurable function $b$ that satisfies 
\begin{equation}\label{excond1}
b(s)\ge\beta(s,y,E)\quad \text{for every}\quad s,y\in [0,m]\quad \text{and}
\quad E\in L^1_+(0,m)
\end{equation} 
and there exists a $R>0$ such that for $||H||_{L^1}+P^0>R$ we have
\begin{equation}\label{excond2}
\int_0^m\kappa(s)F_b(s,H(s),P^0)\,\ud s\le c,
\end{equation}
where
\begin{equation*}
F_b(s,H(s),P^0)=\int_0^s\exp\left\{-\int_x^s\frac{\mu(r,H(r))+\gamma_s(r,P^0)}{
\gamma(r,P^0)}\,\ud
r\right\}\frac{b(x)}{\gamma(x,P^0)}\,\ud x,
\end{equation*}
and $c$ satisfies
\begin{equation}\label{cinequality2}
\kappa(s)\ge c\,b(s),
\end{equation}
for every $s\in [0,m]$. 
Moreover assume that there exists a separable underestimator of the fertility,
\begin{equation}\label{excond3a}
0\le \beta_1(s)\bar{\beta}_1(y,E(y))\le\beta(s,y,E(y)) 
\end{equation} 
for every $s,y\in [0,m],\,E\in L^1_+(0,m)$ such that 
$\bar{\beta}_1$ together with $F_1$ satisfies 
\begin{equation}\label{excond3}
\int_0^m\bar{\beta}_1(s,\mathbf{0})F_1(s,\mathbf{0},0)\,\ud s>1.
\end{equation} 
Then model \eqref{equation1}-\eqref{equation4} admits a positive steady state
$p_*\in W^{1,1}(0,m)$.
\end{theorem}
\noindent {\bf Proof.} We begin by replacing $\beta(s,y,E(y))$ by
\begin{equation*}
\beta(s,y,E(y))-\beta_1(s)\bar{\beta}_1(y,E(y))\ge 0 
\end{equation*}
from condition \eqref{excond3a} and then decompose this remainder as indicated
in equation \eqref{betaseries}. We also note that the decomposition of the
$\beta^l$ functions in 
\eqref{betaseries} can be achieved 
in the way that there is a common first term for every $l$, i.e.~we may write 
$\beta_1(s)\bar{\beta}_1(y,E(y))=\beta^l_1(s)\bar{\beta}^l_1(y,E(y))$ for every
$l$.  This gives a sequence of  finite
rank approximations
\begin{equation*}
\beta^l(s,y,E(y)) = \beta_1(s)\bar{\beta}_1(y,E(y)) +
\sum_{k=2}^{l}\beta^l_k(s)\bar{\beta}^l_k(y,E(y)).
\end{equation*}
Thus, by Theorem \ref{theorem:ex-sep}, every $\Phi^l$ corresponding to a finite
rank approximant $\beta^l$ has a fixed point $(E^l_*,\mathbf{P}^l_*)$ in the
common shell 
\begin{equation*}
\mathcal{S}_{\rho,r}=\{x\in\mathcal{K}\,:\,\rho\le||x||_{\mathcal{X}}\le r\}.
\end{equation*}
Indeed, due to the uniform bound  from inequalities \eqref{excond1} and
\eqref{excond2}, the outer radius $r$ can be chosen uniformly. The common lower
radius $\rho$ can be guaranteed since all approximants $\beta^l$ begin with a
common first term for which condition \eqref{sscond1} from Theorem
\ref{theorem:ex-sep} holds. By the additional gain in regularity due to
$\Phi^l_1$, we have that $E^l_* \in W^{1,1}(0,m)$. By straightforward
calculations, 
\begin{equation*}
\begin{aligned}
& F_j'(s,E^l_*(s),\mathbf{P}^l_*)=\frac{\beta_j(s)}{\gamma(s,P^{0,l}_*)}+
\frac{\mu(s,E^l_*(s))+\gamma_s(s,P^{0,l}_*)}{\gamma(s,P^{0,l}_*)}F_j(s,E^l_*(s),
\mathbf{P}^l_*), \\
& F_j''(s,E^l_*(s),\mathbf{P}^l_*)=\frac{\beta_j'(s)}{\gamma(s,P^{0,l}_*)}-\frac{\beta_j(s)\gamma_s(s,P^{0,l}_*)}{\gamma^2(s , P^{0,l}_*)} \\
&\phantom{=}+\frac{\mu_s(s,E^l_*(s))+\mu_E(s,E^l_*(s))(E^l_*)'(s)+\gamma_{ss}(s ,P^{0,l}_*)}{\gamma(s,P^{0,l}_*)}F_j(s,E^l_*(s),\mathbf{P}^l_*) \\
&\phantom{=}-\frac{\gamma_s(s,P^{0,l}_*)\left(\mu(s,E^l_*(s))+\gamma_s(s,P^{0,l}_*)\right)}{\gamma^2(s,P^{0,l}_*)}F_j(s,E^l_*(s),\mathbf{P}^l_*) \\
&\phantom{=}+\frac{\mu(s,E^l_*(s))+\gamma_s(s,P^{0,l}_*)}{\gamma(s,P^{0,l}_*)}F_j'(s,E^l_*(s),\mathbf{P}^l_*)
\end{aligned}
\end{equation*}
This shows that
$F_j(\,\cdot\,,E^l_*(\,\cdot\,),\mathbf{P}^l_*)\in W^{2,1}(0,m)$ for all $l$ and
$j=1,\dots,l$, and moreover, that this family is uniformly bounded. For every
$l$ the fixed point yields a steady state of the approximate problem by
\begin{equation*}
p_*^l(s)=\sum_{j=1}^lP_*^{j,l}F_j(s,E_*^l(s),P_*^{0,l}).
\end{equation*}
Since $p_*^l$ is a linear
combination of elements in $W^{2,1}(0,m)$, it is itself in $W^{2,1}(0,m)$. 
The uniform bound on the
derivatives in 
\eqref{excond4} implies that
\begin{equation*}
||p_*^l||_{W^{2,1}}=\left|\left|p^l_*\right|\right|_{L^1}+\left|\left|\frac{d\,
p^l_*}{ds}\right|\right|_{L^1}
+\left|\left|\frac{d^2\,p^l_*}{ds^2}\right|\right|_{L^1}\le \tilde{M}\quad
\text{for
all}\quad l\in\mathbb{N}.
\end{equation*} 
This means that $\left\{p_*^l\right\}_{l=1}^{\infty}$ is a bounded set in       
$W^{2,1}(0,m)$ which is compactly embedded in $W^{1,1}(0,m)$ (see e.g.~Theorem
6.2 in \cite{Adams}). Therefore the sequence $p_*^l$ has a convergent
subsequence, again denoted by $p_*^l$, with limit point $p_*$ in
$W^{1,1}(0,m)$. This $p_*$ is the natural candidate for a positive steady state
of model \eqref{equation1}-\eqref{equation4}. 
Next we show that
\begin{equation*}
\lim_{l\to\infty}\left|\left|\Psi(p_*^l)\right|\right|_{L^1}=0,
\end{equation*}
where $\Psi$ is the nonlinear operator defined in \eqref{eq1}. Let $\Psi^l$
denote the nonlinear operator corresponding to the partial sum fertility
function $\beta^l$. Then 
\begin{align*}
& \left|\left|\Psi(p_*^l)-\Psi^l(p_*^l)\right|\right|_{L^1}  \\
& \quad\quad\le\int_0^m\left|\int_0^m\left[\beta(s,y,E^l_*(y))-\beta^l(s,y,
E^l_*(y))\right]p_*^l(y)\,\ud y\right|\,\ud s \\
& \quad\quad\le\int_0^m\left|\left|\left[\beta(s,\,\cdot\,,E^l_*(\,\cdot\,
))-\beta^l(s,\,\cdot\,,E^l_*(\,\cdot\,))\right]p_*^l(\,\cdot\,)
\right|\right|\,\ud s  \\
& \quad\quad\le ||p_*^l||_{L^1}\int_0^m\left|\left|\beta(s,\,\cdot\,,
E^l_*(\,\cdot\,)) -\beta^l(s,\,\cdot\,,E^l_*(\,\cdot\,))
\right|\right|_{L^\infty}\,\ud s \\
& \quad\quad\le K \left|\left|\beta(\,\cdot\,,\,\cdot\,,E_*^l(\,\cdot\,))
-\beta^l(\,\cdot\, ,\,\cdot\,,E_*^l(\,\cdot\,))\right|\right|_{L^\infty},
\end{align*}
for some positive constant $K$. Finally we have
\begin{equation*}
\left|\left|\Psi(p_*)\right|\right|_{L^1}=\left|\left|\Psi\left(\lim_{l\to\infty
}p_*^l\right)\right|\right|_{L^1}=
\left|\left|\lim_{l\to\infty}\Psi\left(p_*^l\right)\right|\right|_{L^1}=\lim_{
l\to\infty}\left|\left|\Psi\left(p_*^l\right)\right|\right|_{L^1}=0,
\end{equation*}
where the second equality follows from the continuity of the operator $\Psi$.
Thus $p_*$ is the desired steady state.
\hfill $\Box$

\begin{remark}
We note that conditions \eqref{excond1}-\eqref{excond3} are natural and
biologically relevant. They are similar to the ones obtained in \cite{P} for the
existence of a positive steady state of a nonlinear age-structured model. For
our model the introduction of a net reproduction function (which will be an
operator) will be somewhat cumbersome and biologically less straightforward, 
see the next section. However, it is still shown that conditions
\eqref{excond1}-\eqref{excond3}  require that the growth rate of the population
is larger than one close to the zero steady state while the growth rate of the
population is small (definitely less than 1) for large population sizes. 
\end{remark}

\section{Asymptotic behavior}\label{section:Asymptotic}

In this section we will investigate the asymptotic behavior of solutions. 
Our approach is based on a formal linearization around a steady state solution
and on a careful spectral analysis of the linearized operator. More precisely,
our goal is to establish conditions which guarantee that the growth bound
$\omega_0$ of the linearized semigroup is negative, respectively positive. 
We will show that the growth bound of the linearized semigroup can be
completely characterized by the spectrum of its generator. However, as we
noted before, the main difficulty is that the eigenvalues of the linearized
semigroup generator and therefore the spectral bound cannot be
cha\-rac\-te\-ri\-zed directly via eigenvalues, in the general case. 

We note that the Principle of Linearized Stability has so far only been 
established for semilinear models, see e.g.~\cite{GuMC,P,WEB}, but not for
general quasilinear equations, such as the one we treat in this paper.
Therefore, for the remainder of the paper, we make the additional assumption that the growth rate 
is separable
\begin{equation}\label{sep_gamma}
\gamma(s,P)=\gamma_1(s)\gamma_2(P). 
\end{equation}
This is plausible from the biological point of view, as the growth rate is
modulated by the total weighted population, equally for individuals of all
sizes. Then the quasilinear problem \eqref{equation1}-\eqref{equation4} can be written in the form
\begin{equation}\label{quasilinear}
\frac{dp}{dt} = g(p)\mathbb{A}p +F(p), \quad p(0) = p_0,
\end{equation}
where 
\begin{equation*}
g(p):=\gamma_2(P)=\gamma_2\left(\int_0^m \kappa(s) p(s,t)\,\ud
s\right), \quad
\mathbb{A}p=\frac{\partial}{\partial s}\left(\gamma_1 p\right),
\end{equation*}
and the recruitment and mortality terms in equation \eqref{equation1} are
incorporated in the nonlinear operator $F$. 
Grabosch and Heijmans in \cite{GH} introduced the transformation
\begin{equation*}
\tau_p(t) = \int_0^t g(p(s))\,\ud s
\end{equation*}
and defined
\begin{equation*}
q(\tau) = p(t_p(\tau)), \quad  \textrm{for}\: \tau \ge0,
\end{equation*}
where $t_p$ is the inverse function of $\tau_p$. It is then verified that $q$ 
satisfies the semilinear equation
\begin{equation}\label{semilinear}
\frac{dq}{d\tau} = \mathbb{A}q(\tau) +B(q(\tau)), \quad q(0) = p_0, 
\end{equation}
with the same initial value as in \eqref{quasilinear}, where the nonlinear operator $B$ is defined via 
$B(q)=F(q)/g(q)$. This requires that $g$ is a continuous, strictly positive and bounded function, as
is guaranteed by our assumptions \eqref{assumptions}. In \cite{GH} Grabosch and
Heijmans showed that
\begin{enumerate}
 \item Solutions of problems \eqref{quasilinear} and \eqref{semilinear} are in
one-to-one correspondence with each other \cite[Theorem 3.4]{GH}.
\item  Problems \eqref{quasilinear} and \eqref{semilinear} have the same
equilibrium solutions, and an equilibrium solution of \eqref{quasilinear} is stable if and only if
it is stable for \eqref{semilinear} \cite[Theorem 5.1]{GH}.
\end{enumerate}
In what follows, we make the assumption of a separable growth rate \eqref{sep_gamma}, but keep the
notation $\gamma(s,P)$ to avoid cumbersome expressions.

\subsection{Linearization around steady states}\label{subsection:Linearization}
Given a stationary (time independent) solution $p_*$ of system
\eqref{equation1}-\eqref{equation4}, we introduce the perturbation $u=u(s,t)$ of
$p$ by 
making the ansatz $p=u+p_*$ and we substitute this into equations
\eqref{equation1}-\eqref{equation4}. 
Then we are using Taylor series expansions of the vital rates of the following
form
\begin{equation*}
f(x,E)=f(x,E_*)+f_E(x,E_*)(E-E_*)+\text{``higher order terms''},
\end{equation*}
to arrive at
\begin{equation}\label{prelinear}
\begin{aligned}
&u_t(s,t) 
+\big(\left(\gamma(s,P_*)+\gamma_P(s,P_*)U(t)\right)\left(u(s,
t)+p_*(s)\right)\big)_s  \\
&\:= 
-\big(\mu(s,E(s,p_*))+\mu_E(s,E(s,p_*))E(s,u)\big)\left(u(s,t)+p_*(s)\right)
\\
&\: +\int_0^m\big(\beta(s,y,E(y,p_*))+\beta_E(s,y,E(y,p_*))E(y,u)\big)
\left(u(y, t)+p_*(y)\right)\,\ud y,
\end{aligned}
\end{equation}
where we have defined 
\begin{equation*}
U(t)=\int_0^m\kappa(s)u(s,t)\,\ud s.
\end{equation*} 
Next we omit the nonlinear terms in equation \eqref{prelinear} to arrive at the
linearized problem  
\begin{equation}\label{linear1}
\begin{aligned}
& u_t(s,t)+\gamma(s,P_*)u_s(s,t)+\gamma_s(s,P_*)u(s,t)+\gamma_{Ps}(s,
P_*)p_*(s)U(t) \\
&\phantom{=} +\gamma_P(s,P_*)p'_*(s)U(t) \\
&=-\mu(s,E(s,p_*))u(s,t)-\mu_E(s,E(s,p_*))p_*(s)E(s,
u) \\
&\phantom{=}
+\int_0^m\bigg(\beta(s,y,E(y,p_*))u(y,t)+\beta_E(s,y,E(y,p_*))p_*(y)E(y,
u)\bigg)\, \ud y,
\end{aligned}
\end{equation}
with the boundary condition
\begin{equation}\label{linear2}
\gamma(0,P_*)u(0,t)=0.
\end{equation}
Equations \eqref{linear1}--\eqref{linear2} are accompanied by the initial
condition
\begin{equation}\label{linear3}
    u(s,0)=u_0(s).
\end{equation}
The linearized problem \eqref{linear1}-\eqref{linear3} is treated effectively in the framework of semigroup theory. 
To this end we cast the linearized system \eqref{linear1}-\eqref{linear3} in the form
of an abstract Cauchy problem on the state space 
$L^1(0,m)$ as follows:
\begin{equation*}
    \frac{d}{dt}\, u = \left({\mathcal A} + {\mathcal B}+ {\mathcal C}+{\mathcal
D}+{\mathcal F}\right)\,u,\quad    u(s,0)=u_0(s),
\end{equation*}
where 
\begin{align}
{\mathcal A} u = & -\gamma(\cdot,P_*)\,u_s\,\,
\text{with}\,\,\text{Dom}({\mathcal A})=\left\{u\in
W^{1,1}(0,m)\,|\,u(s=0)=0\right\}, \label{aop}\\
{\mathcal B} u = &
-\left(\gamma_s(\cdot,P_*)+\mu(\cdot,E(\cdot,p_*))\right)u,\label{bop}\\
{\mathcal C} u = &
-\left(\gamma_{Ps}(\cdot,P_*)p_*(\cdot)+\gamma_P(\cdot,
P_*)p'_*(\cdot)\right)\int_0^m\kappa(y)u(y)\,\ud y \nonumber \\
= & -\rho_*(\cdot)\int_0^m\kappa(y)u(y)\,\ud y,\label{cop} \\
{\mathcal D} u= & -\mu_E(\cdot,E(\cdot,p_*))p_*(\cdot)E(\cdot,u),\label{dop} \\
{\mathcal F} u= & \int_0^m\beta(\cdot,y,E(y,p_*))u(y)\,\ud
y+\int_0^m\beta_E(\cdot,y,E(y,p_*))p_*(y)E(y,u)\,\ud y.\label{fop}
\end{align}
\begin{proposition}\label{existencelinear}
The operator $\mathcal{A+B+C+D+F}$ generates a strong\-ly continuous
quasicontractive semigroup 
$\left\{\mathcal{T}(t)\right\}_{t\ge 0}$ of bounded linear operators on
$L^1(0,m)$, which is positive if the operator $\mathcal{C+D+F}$ is positive.
\end{proposition} 
We omit the proof of the above  proposition as it can be established following
similar results in \cite{FGH,FH2,FH3,FHin}.

\subsection{Stability results for positive
equilibria}\label{subsection:Stability}
In this section we establish linear sta\-bi\-li\-ty/instability results for a
general steady state $p_*$. We will treat the extinction steady
state in the next subsection in detail. When addressing stability questions 
the main difficulty is that, in general, the point spectrum of the semigroup generator cannot be
characterized via zeros of an associated characteristic function.
We will overcome this problem by using positive perturbation arguments.

\begin{theorem}\label{instability}
Assume that 
\begin{equation}\label{instabcond1}
\int_0^m\kappa(s)\int_0^s\exp\left\{-\int_y^s\frac{\gamma_s(x,P_*)+\mu(x,E(x,
p_*))}{\gamma(x,P_*)}\,\ud x\right\}
\frac{\rho_*(y)}{\gamma(y,P_*)}\,\ud y\,\ud s<-1.
\end{equation}
Furthermore assume that
\begin{align}
& \rho_*(s)\le 0,\quad s\in [0,m]\quad \text{and}\quad
\exists\,\varepsilon>0\,\, \text{such that}\,\, \rho_*(s)\ne 0\, \text{for
a.e.}\, s\in[0,\varepsilon], \label{instabcond2} \\
& \mu_E(s,E(s,p_*))\le 0,\quad \beta_E(s,y,E(y,p_*))\ge 0,\quad y,s\in
[0,m].\label{instabcond3}
\end{align}
Then the steady state $p_*$ of model \eqref{equation1}-\eqref{equation4} is
linearly unstable.
\end{theorem}
\noindent {\bf Proof.} 
First we note that $\mathcal{A+B+C}$ generates a positive, irreducible and
eventually compact semigroup if conditions \eqref{instabcond2} 
hold true. Eventual compactness is easily shown, see Proposition \ref{reg}. 
To establish irreducibility we consider the resolvent equation 
\begin{equation*}
R(\lambda,\mathcal{A+B+C})h=u,\quad h\in\mathcal{X}^+.   
\end{equation*}
The solution is obtained as:
\begin{align}
u(s)= &
\int_0^s\exp\left\{-\int_y^s\frac{\gamma_s(x,P_*)+\mu(x,E(x,p_*))+\lambda}{
\gamma(x,
P_*)}\,\ud x\right\}\frac{h(y)}{\gamma(y,P_*)}\,\ud y\nonumber \\
&
-U\int_0^s\exp\left\{-\int_y^s\frac{\gamma_s(x,P_*)+\mu(x,E(x,p_*))+\lambda}{
\gamma(x,
P_*)}\,\ud x\right\}\frac{\rho_*(y)}{\gamma(y,P_*)}\,\ud y,\label{irredeq1}
\end{align}
where 
\begin{align}
U= & \int_0^m\kappa(s)u(s)\,\ud s \nonumber \\
= &
\frac{\int_0^m\kappa(s)\int_0^s\exp\left\{-\int_y^s\frac{\gamma_s(x,P_*)+\mu(x,
E(x,p_*))+\lambda}{\gamma(x,
P_*)}\,\ud x\right\}\frac{h(y)}{\gamma(y,P_*)}\,\ud y\,\ud
s}{1+\int_0^m\kappa(s)\int_0^s\exp\left\{-\int_y^s\frac{\gamma_s(x,P_*)+\mu(x,
E(x,p_*))+\lambda}{\gamma(x,P_*)}\,\ud
x\right\}\frac{\rho_*(y)}{\gamma(y,P_*)}\,\ud y\,\ud s}.
\end{align}
The second equality is obtained from \eqref{irredeq1} by multiplying by $\kappa$
and integrating from $0$ to $m$. Hence for $\lambda$ large enough 
$U>0$ and $u\gg 0$ follows from condition \eqref{instabcond2}. The
irreducibility and eventual compactness of the semigroup imply that the spectrum
$\sigma(\mathcal{A+B+C})$ is not empty, see e.g. Th. 3.7 in Sect. C-III in
\cite{AGG}, hence the spectrum $\sigma(\mathcal{A+B+C})$ contains at least one 
eigenvalue. Next we find the solution of the eigenvalue equation
\begin{equation*}
(\A+\B+\mathcal{C}) u = \lambda u
\end{equation*}
as
\begin{equation}\label{eigv}
u(s)
=-U\int_0^s\exp\left\{-\int_y^s\frac{\gamma_s(x,P_*)+\mu(x,E(x,p_*))+\lambda}{
\gamma(x,
P_*)}\,\ud x\right\}\frac{\rho_*(y)}{\gamma(y,P_*)}\,\ud y.
\end{equation}
Then we multiply the solution \eqref{eigv} by $\kappa$ and integrate over
$[0,m]$ to obtain
\begin{equation*}
\begin{aligned}
&U= -U  \int_0^m \kappa(s) \\ 
&\int_0^s\exp\left\{-\int_y^s\frac{\gamma_s(x,P_*)+\mu(x,E(x,
p_*))+\lambda}{\gamma(x,P_*)}\,\ud x\right\}\frac{\rho_*(y)}{\gamma(y,P_*)}\,\ud
y\,\ud s.
\end{aligned}
\end{equation*}
We note that, if $U=0$ then equation \eqref{eigv} shows that $u(s)\equiv 0$,
hence we have a non-trivial eigenvector if and only if 
$U\ne 0$ and $\lambda$ satisfies the following characteristic equation
\begin{equation}\label{char1}
\begin{aligned}
&1 =K(\lambda)  \eqdef \\
&\int_0^m \kappa(s)\int_0^s\exp\left\{-\int_y^s\frac{\gamma_s(x,P_*)+\mu(x,E(x,
p_*))+\lambda}{\gamma(x,P_*)}\,\ud x\right\}\frac{\rho_*(y)}{\gamma(y,P_*)}\,\ud
y\,\ud s.
\end{aligned}
\end{equation}
It is easily shown that
\begin{equation*}
\lim_{\lambda\ra +\infty} K(\lambda)  = 0,
\end{equation*}
therefore it follows from the assumption $K(0)>1$ \eqref{instabcond1}, on the
grounds of the Intermediate Value Theorem, that
equation \eqref{char1} has a positive (real) solution. 
Hence we have
\begin{equation*}
0<s(\mathcal{A+B+C}),
\end{equation*}
where $s(\mathcal{A+B+C})$ stands for the spectral bound of the operator
$\mathcal{A+B+C}$. The operators 
$\mathcal{D}$ and $\mathcal{F}$ are positive if conditions \eqref{instabcond3}
hold true. 
We have for the spectral bound (see e.g.~Corollary VI.1.11 in \cite{NAG})
\begin{equation*}
0<s(\mathcal{A+B+C})\le s(\mathcal{A+B+C+D+F}).
\end{equation*}
Since the growth bound of the semigroup is bounded below by the spectral bound
of its generator, the proof is completed.
\hfill $\Box$

We note that the instability conditions \eqref{instabcond3} imply that mortality
is a non-increasing function of the environment and fertility is a
non-decreasing function of the environment.  

Next we establish conditions which guarantee that the equilibrium solution $p_*$
is linearly asymptotically stable. To this end we establish first that the
spectrum of the semigroup generator $\mathcal{A+B+C+D+F}$ consists of
eigenvalues only and that the spectral mapping theorem holds true. Then the
growth bound of the semigroup equals the spectral bound of its ge\-ne\-ra\-tor.
These follow however from the following result (see e.g.~Corollary IV.3.12 in
\cite{NAG}).
\begin{proposition}\label{reg}
\begin{sloppypar}
The semigroup $\{{\mathcal T}(t)\}_{t\geq 0}$ generated by the operator
\mbox{${\mathcal A} + {\mathcal B}+ {\mathcal C}+ {\mathcal D}+\mathcal{F}$} 
is eventually compact.
\end{sloppypar}
\end{proposition}
\noindent {\bf Proof.}
We only sketch the proof here since analogous results for simpler problems can
be found in \cite{FGH, FH2,FH3,FHin}. Due to the zero flux boundary condition
and the finite maximal size, the operator $\mathcal{A+B}$ generates a nilpotent
semigroup. The biological interpretation is that in the absence of recruitment
the population dies out independently of the initial condition. $\mathcal{C}$ is
a bounded linear operator of rank one, hence it is compact. It only remains to
establish that the
bounded linear integral operators $\mathcal{D}$ and $\mathcal{F}$ are compact. 
These however, can be deduced using the Fr\'{e}chet-Kolmogorov compactness
criterion (see e.g.~Chapter X in \cite{Y}) from the regularity assumptions we
made on the model ingredients, see the proof of Theorem 12 in \cite{FGH} for
more details. \hfill $\Box$

The previous result guarantees that stability is determined by the leading
eigenvalue of the semigroup generator 
$\mathcal{A+B+C+D+F}$, unless the spectrum is empty. In that case one further
needs to establish
positivity of the semigroup which guarantees that the growth bound coincides
with the spectral bound, which by definition equals minus infinity. However, as
we noted
before, the eigenvalue equation
\begin{equation*}
\left(\mathcal{A+B+C+D+F}-\lambda\mathcal{I}\right)u=0,\quad\quad u(s=0)=0,
\end{equation*}
cannot be solved explicitly. This is due to the
infinite dimensional non\-li\-ne\-arity in the original problem and to the 
very general recruitment term. 

Compact perturbations do not change the essential spectrum 
of the semigroup generator. Therefore our approach to establish stability is
to find a positive compact perturbation of the 
generator for which we can characterize the point spectrum via zeros of an
associated characteristic function. To this end, 
for a separable fertility function $\beta$, we introduce the following
operators on $L^1(0,m)$,
\begin{align*}
&
\overline{\mathcal{D}}u=-\mu_E(\cdot,E(\cdot,p_*))p_*(\cdot)\int_0^mw(r)u(r)\,
\ud r,  \\
&
\overline{\mathcal{F}}_2u=\beta_1(\cdot)\int_0^m\beta_{2_E}(y,E(y,p_*))p_*(y)\,
\ud y\int_0^mw(r)u(r)\,\ud r. 
\end{align*}
We formulate our main stability result in the following theorem.
\begin{theorem}\label{stability}
Assume that there exists a function 
\begin{equation*}
\widetilde{\beta}(s,y,E(y,p_*))=\beta_1(s)\beta_2(y,E(y,p_*))
\end{equation*} 
such that 
\begin{equation}\label{stabcond3}
\beta(s,y,E(y,p_*))\le\widetilde{\beta}(s,y,E(y,p_*)),\quad s,y\in [0,m].
\end{equation}
Furthermore assume that the following conditions hold true
\begin{align}
& \mu_E(s,E(s,p_*))\le 0,\quad\widetilde{\beta}_E(s,y,E(y,p_*)))\ge 0,\quad
s,y\in [0,m], \label{stabcond1} \\
& \gamma_{Ps}(s,P_*)p_*(s)+\gamma_P(s,P_*)p'_*(s)\le 0, \quad s\in [0,m],
\label{stabcond2}
\end{align} 
and the characteristic function $K_{\widetilde{\beta}}(\lambda)$ 
given by equation \eqref{def_K} corresponding to $\widetilde{\beta}$ and to the
modified 
operators $\overline{\mathcal{D}}$ and $\overline{\mathcal{F}}_2$ does not have
a zero with non-negative real part. 
Then the stationary solution $p_*$ is linearly asymptotically stable.
\end{theorem}
\noindent {\bf Proof.} 
We obtain the solution of the eigenvalue equation
\begin{equation}\label{eigveq}
\left(\mathcal{A+B+C}+\overline{\mathcal{D}}+\mathcal{F}_1+\overline{\mathcal{F}
}_2-\lambda\mathcal{I}\right)u=0,
\end{equation}
as
\begin{equation}\label{soleigveq}
\begin{aligned}
u(s)= &  U_1\int_0^s f_0^\lambda(s,y)f_1(y)\,\ud y  \\
& +U_2\int_0^sf_0^\lambda(s,y)f_2(y)\,\ud y+U_3\int_0^s
f_0^\lambda(s,y)f_3(y)\,\ud y,
\end{aligned}
\end{equation}
where
\begin{align*}
& U_1=\int_0^m\kappa(s)u(s)\,\ud s,\quad U_2=\int_0^mw(s)u(s)\,\ud
s,\\
& U_3=\int_0^m\beta_2(s,E(s,p_*))u(s)\,\ud s, \\
&
f_0^\lambda(s,y)=\exp\left\{-\int_y^s\frac{\gamma_s(x,P_*)+\mu(x,E(x,
p_*))+\lambda}{
\gamma(x,P_*)}\,\ud x\right\}, \\
& f_1(y)=\frac{-\gamma_{Ps}(y,P_*)p_*(y)-\gamma_P(y,P_*)p'_*(y)}{\gamma(y,P_*)},
\quad f_3(y)=\frac{\beta_1(y)}{\gamma(y,P_*)},\\
& f_2(y)=\frac{\beta_1(y)\int_0^m\beta_2(x,E(x,p_*))p_*(x)\,\ud
x-\mu_E(y,E(y,p_*))p_*(y)}{\gamma(y,P_*)}.
\end{align*}
We multiply equation \eqref{soleigveq} by $\kappa,w$ and by $\beta_2$ and
integrate from $0$ to $m$, respectively to obtain
\begin{align}
& U_1(1+a_{11}(\lambda))+U_2a_{12}(\lambda)+U_3a_{33}(\lambda)=0,\label{homog1}
\\
& U_1a_{21}(\lambda)+U_2(a_{22}(\lambda)+1)+U_3a_{23}(\lambda)=0,\label{homog2}
\\
& U_1a_{31}(\lambda)+U_2a_{32}(\lambda)+U_3(a_{33}(\lambda)+1)=0, \label{homog3}
\end{align}
where for $i=1,2,3$,
\begin{align*}
& a_{1i}(\lambda)=\int_0^m\kappa(s)\int_0^sf_0^\lambda(s,y)f_i(y)\,\ud y\,\ud
s,  \\
& a_{2i}(\lambda)=\int_0^mw(s)\int_0^sf_0^\lambda(s,y)f_i(y)\,\ud y\,\ud s,   \\
& a_{3i}(\lambda)=\int_0^m\beta_2(s,E(s,p_*))\int_0^sf_0^\lambda(s,y)f_i(y)\,\ud
y\,\ud s. 
\end{align*}
If
$\lambda\in\sigma(\mathcal{A+B+C}+\overline{\mathcal{D}}+\mathcal{F}_1+\overline
{\mathcal{F}_2})$ 
then the eigenvalue equation \eqref{eigveq} admits a non-trivial solution $u$
hence there exists a non-zero vector 
$(U_1,U_2,U_3)$ which solves equations \eqref{homog1}-\eqref{homog3}. 
To the contrary, if $(U_1,U_2,U_3)$ is a non-zero solution of
equations \eqref{homog1}-\eqref{homog3} for some 
$\lambda\in\mathbb{C}$ then \eqref{soleigveq} yields a non-trivial solution
$u$. 
This is because the only scenario for $u$ to vanish would yield
\begin{align*}
-U_1\int_0^sf_0^\lambda(s,y)f_1(y)\,\ud y
&=U_2\int_0^sf_0^\lambda(s,y)f_2(y)\,\ud y \\
&\phantom{=} +U_3\int_0^sf_0^\lambda(s,y)f_3(y)\,\ud y,
\end{align*}
for every $s\in [0,m]$. This however, together with conditions
\eqref{stabcond1}-\eqref{stabcond2} would yield $U_1=U_2=U_3=0$, a
contradiction.
Thus $\lambda\in\mathbb{C}$ is an eigenvalue of
$\mathcal{A+B+C}+\overline{\mathcal{D}}+\mathcal{F}_1+\overline{\mathcal{F}_2}$ 
if and only if it satisfies the characteristic equation
\begin{equation}\label{def_K}
K_{\widetilde{\beta}}\,(\lambda)\eqdef\det
\left( \begin{array}{llll}
1+a_{11}(\lambda) & a_{12}(\lambda) & a_{13}(\lambda)\\
a_{21}(\lambda) & 1+a_{22}(\lambda) & a_{23}(\lambda) \\
a_{31}(\lambda) & a_{32}(\lambda) & 1+a_{33}(\lambda)
\end{array}\right)=0.
\end{equation}
Next we observe that conditions \eqref{stabcond1}-\eqref{stabcond2} guarantee
that both $\mathcal{C}$, $\mathcal{D}$ and 
$\mathcal{F}$ are positive operators. Therefore we conclude that
$\mathcal{A+B+C+D+F}$ is a generator of a positive semigroup. 
Moreover, the operators $(\overline{\mathcal{D}}-\mathcal{D})$, 
$(\mathcal{F}^{\widetilde{\beta}}_1-\mathcal{F}_1^\beta)$ and 
$(\overline{\mathcal{F}}_2^{\widetilde{\beta}}-\mathcal{F}_2^\beta)$ are all
positive and bounded. 
We have
\begin{align*}
& s(\mathcal{A+B+C+D+F})
=s\left(\mathcal{A+B+C+D}+\mathcal{F}_1^{\beta}+\mathcal{F}_2^\beta\right) \\
& \quad\le 
s\left(\mathcal{A+B+C+D}+\mathcal{F}_1^{\beta}+\mathcal{F}_2^\beta+\mathcal{
\overline{D}}-\mathcal{D}+\mathcal{F}^{\widetilde{\beta}}_1-\mathcal{F}
_1^\beta+\overline{\mathcal{F}}_2^{\widetilde{\beta}}-\mathcal{F}_2^\beta\right)
 \\
& \quad
=s\left(\mathcal{A+B+C}+\overline{\mathcal{D}}+\mathcal{F}_1^{\widetilde{\beta}}
+\overline{\mathcal{F}}_2^{\widetilde{\beta}}\right)<0,
\end{align*} 
and the proof is completed.
\hfill $\Box$

\begin{example}
The crucial assumption of the previous theorem is that the characteristic
function $K_{\widetilde{\beta}}$ does not have a root 
with non-negative real part. Here we only present an example when this condition
may be easily verified. In particular, let us assume that
\begin{equation*}
\kappa\equiv c_1,\quad w\equiv c_2,\quad \beta_2\equiv c_3
\end{equation*} 
for some constants $c_1,c_2,c_3>0$. In this special case the characteristic
equation $K_{\widetilde{\beta}}(\lambda)=0$ takes the simple form
\begin{equation*}
1=\int_0^m\int_0^sf_0^\lambda(s,y)\big(c_1f_1(y)+c_2f_2(y)+c_3f_3(y)\big)
\ud y\,\ud s.
\end{equation*}
Therefore there exists a unique dominant real eigenvalue, which is negative, if
\begin{align*}
&
\int_0^m\int_0^s\exp\left\{-\int_y^s\frac{\gamma_s(x,P_*)+\mu(x,E(x,p_*))}{
\gamma(x,P_*)}\,\ud x\right\}  \\
&
\quad\quad\times\left(\frac{\beta_1(y)\left(c_3+c_2\int_0^m\beta_2(s,E(s,
p_*))p_*(s)\,\ud s\right)}{\gamma(y,P_*)}\right.  \\
&
\quad\quad\quad\quad\left.-\frac{c_2\mu_E(y,E(y,p_*))+c_1\gamma_{Ps}(y,
P_*)p_*(y)+c_1\gamma_P(y,P_*)}{\gamma(y,P_*)}\right)\,\ud y\,\ud s \\ 
& \quad <1.
\end{align*}
The latter condition may be easily verified for fixed model ingredients.
\end{example}

\subsection{The extinction equilibrium}\label{subsection:Extinction}
In this subsection we establish a simple criterion for the
stability/instability of the trivial steady state $p_*\equiv 0$. The stability
of the trivial steady state is important from the biological point of view, as
it can answer the question for example if a species can be introduced
successfully into (or can invade) a new habitat. In case of the trivial steady
state the eigenvalue problem can be written as
\begin{equation}\label{simpleeigv}
\left(\mathcal{A+B+F}_{\beta}-\lambda\mathcal{I}\right)u=0,\quad u(s=0)=0,
\end{equation}
where $\lambda\in\mathbb{C}$, the operators $\A$ and $\B$ are defined
via \eqref{aop}-\eqref{bop} with $p_*\equiv 0$ (and $P_*=0$) and
\begin{equation*}
\mathcal{F}_\beta u=\int_0^m\beta(\,\cdot\,,y,\mathbf{0})u(y)\,\ud y.
\end{equation*}
We recall that in case of simpler scramble competition models the so called
net reproduction function played a crucial role in the stability analysis of
equilibrium solutions, see e.g.~\cite{FH}. In particular, we managed to relate
our stability results to a biologically meaningful model ingredient, namely the
net reproduction rate. In case of hierarchical contest competition models this
is less straightforward as we have seen. Nevertheless at least for a separable
fertility function
\begin{equation*}
\beta(s,y,E(y,\cdot))=\beta_1(s)\beta_2(y,E(y,\cdot)),
\end{equation*}  
we may define a {\it net reproduction functional} 
\begin{equation*}
R:L^1_+(0,m)\to \mathbb{R}
\end{equation*}  
of the standing population $p$ as
\begin{align*}
&R(p)  \\ &=
\int_0^m\int_0^s\exp\left\{-\int_y^s\frac{\gamma_s(x,P)+\mu(x,E(x,p))}{
\gamma(x,P)}\,\ud x\right\}\frac{\beta_1(y)\beta_2(s,E(s,p))}{\gamma(y,P)}\,\ud
y\,\ud s, \\
&
=\int_0^m\frac{\beta_2(s,E(s,p))}{\gamma(s,P)}\int_0^s\beta_1(y)\exp\left\{
-\int_y^s\frac{\mu(x,E(x,p))}{
\gamma(x,P)}\,\ud x\right\}\,\ud y\,\ud s, 
\end{align*}
where $P$ is the weighted population according to \eqref{equation4}.
The value of the functional $R$ is the expected number of
offspring to be produced by an individual in her lifetime. Individuals of size
$y$ are produced at a rate $\beta_1$ and need to survive from size $y$ to
size $s$ to reproduce.
We also note that $R(p_*)=1$ is a necessary condition for the existence of a
positive steady
state $p_*$ of our model. It is however not a sufficient condition in contrast
to scramble competition models, i.e.~when density dependence is incorporated via
finite dimensional nonlinearities in the vital rates.

\begin{proposition}\label{trivialss}
Assume that there exists a function $\beta^l$ such that 
\begin{equation*}
\beta^l(s,y)=\beta^l_1(s)\beta^l_2(y)\le\beta(s,y,\mathbf{0})
\end{equation*}
and $R_{\beta^l}(\mathbf{0})>1$ (where $R_{\beta^l}$ stands for the net
reproduction
functional corresponding to the fertility function $\beta^l$). Then the trivial
steady state is linearly unstable. On the other hand if there exists a function
$\beta^u$ with
\begin{equation*}
\beta^u(s,y)=\beta^u_1(s)\beta^u_2(y)\ge\beta(s,y,\mathbf{0})
\end{equation*}
and $R_{\beta^u}(\mathbf{0})<1$, then the trivial steady state is linearly
stable.
\end{proposition}
\noindent {\bf Proof.}
The eigenvalue equation \eqref{simpleeigv} has a nontrivial solution if and only
if $\lambda$ satisfies the characteristic equation
\begin{equation}\label{char2}
\begin{aligned}
&1 =K^l(\lambda)  \eqdef \\
&\int_0^m \beta_2^l(s) \\
&
\hspace{4mm}\times\int_0^s\exp\left\{-\int_y^s\frac{\gamma_s(x,P_*)+\mu(x,E(x,
p_*))+\lambda}{\gamma(x,P_*)}\,\ud
x\right\}\frac{\beta_1^l(y)}{\gamma(y,P_*)}\,\ud y\,\ud s.
\end{aligned}
\end{equation}
It can be shown that $\displaystyle\lim_{\lambda\to
+\infty} K^l(\lambda)=0$. Therefore, if $R_{\beta^l}(\mathbf{0})>1$ holds true,
then the characteristic equation \eqref{char2} admits a positive solution
$\lambda$. 
Since $\mathcal{A+B+F}_{\beta^l}$ generates a positive semigroup, and the
operator $\mathcal{F}_\beta-\mathcal{F}_{\beta^l}$ is positive we have 
$0<s(\mathcal{A+B+F}_{\beta^l})\le s(\mathcal{A+B+F}_{\beta})$ and the
instability part follows. 
Similarly, if $R_{\beta^u}(\mathbf{0})<1$ holds true, then the characteristic
equation $1=K^u(\lambda)$ does not have a solution with non-negative real part.
Therefore on the grounds of Proposition \ref{reg}. we have 
\begin{equation*}
s(\mathcal{A+B+F}_\beta+\mathcal{F}_{\beta^u}-\mathcal{F}_\beta)\le
s(\mathcal{A+B+F}_{\beta^u})<0,
\end{equation*} 
and the stability part follows. 
\hfill $\Box$

\section{Concluding remarks}\label{section:Conclusion}

Hierarchical structured population models are important from the application
point of view 
as they describe the evolution of populations in which individuals are
experiencing size-specific environment. 
This is the case for example in a tree population or in a cannibalistic
population. 
They pose a greater mathematical challenge though than scramble competition
models. 
In the recent papers \cite{FH2,FH3} we started to investigate the asymptotic
behavior  
of hierarchical structured partial dif\-fe\-ren\-ti\-al equation models with one state
at birth. 
Nevertheless, the question of the existence of positive equilibrium solutions
for hierarchical models 
remained an open question up to our knowledge. In this work we formulated
biologically relevant 
conditions for the existence of positive steady states of a very general model. 
 
In the second part we focused on the stability of equilibria. As we have also
seen earlier in \cite{FH2,FH3}, the main difficulty to establish easily
verifiable stability/instability conditions for positive steady states 
of hierarchical structured population models is that the point spectrum of the
semigroup generator cannot be characterized explicitly via zeros of an
associated characteristic function. In \cite{FH2, FH3} we devised a
dissipativity calculation to show that the linearized semigroup has negative
growth bound, and therefore the steady state is stable. The stability conditions
we obtained in that way, though, were extremely restrictive. Therefore in this
work we devised a new approach to establish stability by using results from the
theory of compact and positive operators. 

The separability assumption made in \eqref{sep_gamma} is only
required for the result of Grabosch and Heijmans \cite{GH} to be applicable,
and not to prove any of the results concerning the linearized semigroup $\T(t)$.
A Principle of Linearized Stability for fully nonlinear quasilinear
hyperbolic equations remains a problem for future investigations (see also 
\cite{K,TZ,WEB} for related results).

\subsection*{Acknowledgments}

\noindent J.~Z.~Farkas was supported by a personal research grant by the
Carnegie Trust for the Universities of Scotland. Part of this work was done
while J.~Z.~Farkas visited the University of Wisconsin - Milwaukee. Financial
support from the Department of Mathematical Sciences is greatly appreciated.  
P.~Hinow thanks the organizers of the December 2010 Workshop on PDE Models of
Biological Processes at the National Center for Theoretical Sciences at the
National Tsing Hua University in Hsinchu, Taiwan.


\begin{thebibliography}{99}
\bibitem{ADH} (MR2101381)
A.~S. Ackleh, K. Deng, and S. Hu.
\newblock A quasilinear hierarchical size-structured model: well-posedness and
approximation, \newblock {\em Appl. Math. Optim.}\ {\bf 51} (2005), 35--59.

\bibitem{AI} (MR2168831)
A.~S. Ackleh and K. Ito.
\newblock Measure-valued solutions for a hierarchically size-structured
population, \newblock {\em J. Differential Equations}\ {\bf 217} (2005),
431--455.

\bibitem{Adams} (MR2424078)
R. A. Adams, J. J. F. Fournier.
\newblock {\em Sobolev Spaces, 2nd edition}, \newblock Elsevier/Academic Press,
Amsterdam, Boston (2003).

\bibitem{Amann} (MR0415432)
H. Amann. \newblock Fixed point equations and nonlinear eigenvalue problems in
ordered Banach spaces, \newblock {\em SIAM Rev.}\ {\bf 18} (1976), 620--709.

\bibitem{AGG} (MR0839450)
W. Arendt, A. Grabosch, G. Greiner, U. Groh, H.~P. Lotz, U. Moustakas, R. Nagel,
F. Neubrander and U. Schlotterbeck,
\newblock {\em One-Parameter Semigroups of Positive Operators}, \newblock
Springer-Verlag, Berlin, (1986).

\bibitem{CS3} (MR1478297)
\`A. Calsina and J. Salda\~na. \newblock Asymptotic behavior of a model of
hierarchically structured population dynamics, \newblock {\em J. Math. Biol.}\ 
{\bf 35} (1997), 967--987.

\bibitem{CS2} (MR2264557)
\`A. Calsina and J. Salda\~na. \newblock Basic theory for a class of models of
hierarchically structured population dynamics with distributed states in the
recruitment, \newblock {\em Math. Models Methods Appl. Sci.}\  {\bf 16} (2006),
1695--1722.

\bibitem{C1} (MR1293670)
J.~M. Cushing. \newblock The dynamics of hierarchical age-structured
populations, \newblock  {\em J. Math. Biol.}\ {\bf 32} (1994), 705--729.

\bibitem{CUS} (MR1636703)
J.~M. Cushing. \newblock {\em An Introduction to Structured Population
Dynamics}, \newblock SIAM, Philadelphia (1998).

\bibitem{Deimling} (MR0787404)
K. Deimling. \newblock {\em Nonlinear Functional Analysis}, \newblock
Springer-Verlag, Berlin (1985).

\bibitem{NAG} (MR1721989)
K.-J. Engel and R. Nagel. \newblock {\em One-Parameter Semigroups for Linear
Evolution Equations}, \newblock Springer, New York 2000.

\bibitem{FGH} (MR2663320)
J.~Z. Farkas, D.~M. Green and P. Hinow. \newblock Semigroup analysis of
structured
parasite populations, \newblock {\em Math. Model. Nat. Phenom.}\ {\bf 5}
(2010), 94--114.

\bibitem{FH} (MR2285538)
J.~Z. Farkas and T. Hagen. \newblock Stability and regularity results for a
size-structured population
model, \newblock {\em J. Math. Anal. Appl.}\ {\bf 328} (2007), 119--136.

\bibitem{FH2} (MR2552152)
J.~Z. Farkas and T. Hagen. \newblock Asymptotic analysis of a size-structured
cannibalism model with infinite dimensional environmental feedback,
\newblock {\em Commun. Pure Appl. Anal.}\ {\bf 8} (2009), 1825--1839.

\bibitem{FH3} (MR2767889)
J.~Z. Farkas and T. Hagen. \newblock Hierarchical size-structured populations:
The linearized semigroup approach, 
\newblock {\em  Dyn. Contin. Discrete Impuls. Syst. Ser. A Math. Anal.} {\bf 17}
(2010), 639--657.

\bibitem{FHin} (MR2680511)
J.~Z. Farkas and P. Hinow. \newblock On a size-structured two-phase population
model with infinite states-at-birth, \newblock {\em Positivity}\ {\bf 14}
(2010), 501--514.

\bibitem{GH} (MR1076297)
A. Grabosch and H.~J.~A.~M. Heijmans. \newblock Cauchy problems with
state-dependent time evolution, \newblock {\em Japan J. Appl. Math.}\ {\bf 7}
(1990), 433--457.

\bibitem{GuMC} (MR0354068)
M. E. Gurtin and R. C. MacCamy. \newblock Non-linear age-dependent population
dynamics, \newblock {\em Arch. Rational Mech. Anal.} \textbf{54} (1974),
281--300.

\bibitem{HC}
S.~M.~Henson and J.~M. Cushing. \newblock Hierarchical models of intra-specific
competition: scramble versus contest, 
\newblock  {\em J. Math. Biol.}\ {\bf 34} (1996), 755--772.

\bibitem{JC} (MR1972195)
S.~R.-J. Jang and J.~M. Cushing. \newblock A discrete hierarchical model of
intra-specific
competition, \newblock  {\em J. Math. Anal. Appl.}\ {\bf 280} (2003), 102--122.

\bibitem{JC2} (MR2114319)
S.~R.-J. Jang and J.~M. Cushing. \newblock Dynamics of hierarchical models in
discrete-time, 
\newblock  {\em J. Difference Equ. Appl.}\ {\bf 11} (2005), 95--115.

\bibitem{K} (MR1290722)
N. Kato. \newblock A principle of linearized stability for nonlinear evolution
equations, \newblock {\em Trans. Amer. Math. Soc.}\ {\bf 347} (1995),
2851--2868. 

\bibitem{MD} (MR0860963)
J.~A.~J. Metz and O. Diekmann. {\em The Dynamics of Physiologically Structured
Populations}, \newblock  Springer, Berlin, 1986.

\bibitem{P} (MR0702651)
J. Pr\"{u}ss. \newblock On the qualitative behavior of populations with
age-specific interactions, \newblock {\em Comput. Math. Appl.}\ {\bf 9} (1983),
327--339.

\bibitem{TZ} (MR0941101)
S.~L.~Tucker and S.~O.~Zimmerman. \newblock A nonlinear model of population
dynamics containing an arbitrary number of continuous structure variables,
\newblock {\em SIAM J. Appl. Math.}\ {\bf 48} (1988), 549--591.

\bibitem{W2} (MR2593606)
Ch. Walker. \newblock Global bifurcation of positive equilibria in nonlinear
population models, \newblock {\em J. Diff. Eq.}\ {\bf 248} (2010), 1756--1776.

\bibitem{WEB} (MR0772205)
G.~F. Webb. \newblock {\em Theory of Nonlinear Age-Dependent Population
Dynamics}, \newblock Marcel Dekker, New York, (1985).

\bibitem{Y} (MR1336382)
K. Yosida. \newblock{\em Functional Analysis}, \newblock Springer, Berlin,
(1995).

\end{thebibliography}
\end{document}